\documentclass[11pt,letterpaper]{amsart}
\usepackage{amssymb}
\def\AA{{\mathbb A}}

\def\CC{{\mathbb C}}

\def\GG{{\mathbb G}}

\def\PP{{\mathbb P}}
\def\QQ{{\mathbb Q}}

\def\spec{{\rm Spec}}

\def\sing{{\rm sg}}

\def\bs{\backslash}

\def\pt{{\bullet}}
\def\eps{\epsilon}

\def\Acal{{\mathcal A}}

\def\Ccal{{\mathcal C}}

\def\Hcal{{\mathcal H}}

\def\Mcal{{\mathcal M}}

\def\Rcal{{\mathcal R}}
\def\Scal{{\mathcal S}}
\def\Tcal{{\mathcal T}}

\newcommand\pic{\operatorname{Pic}}

\newcommand\proj{\operatorname{Proj}}

\newcommand\sym{\operatorname{Sym}}

\newcommand\GL{\operatorname{GL}}

\newcommand\SL{\operatorname{SL}}

\newtheorem{theorem}{Theorem}[section]
\newtheorem{lemma}[theorem]{Lemma}
\newtheorem{proposition}[theorem]{Proposition}

\newtheorem{corollary}[theorem]{Corollary}

\theoremstyle{definition}

\theoremstyle{remark}

\newtheorem{question}[theorem]{Question}

\title{A perfect stratification of  $\Mcal_g$ for $g\le 5$}
\author{Claudio Fontanari}
\address{Dipartimento di Matematica,
Politecnico di Torino,
Corso Duca degli Abruzzi 24,
10129 Torino (Italia)}
\email{claudio.fontanari@polito.it}
\author{Eduard Looijenga}
\address{Mathematisch Instituut,
Universiteit Utrecht,
P.O.~Box 80.010, NL-3508 TA Utrecht
(Nederland)}
\email{looijeng@math.uu.nl}


\thanks{This project started during a stay of both authors
at the \emph{Institut Mittag--Leffler} (Djursholm, Sverge). 
Support by this institution is gratefully acknowledged.}
\begin{document}

\begin{abstract}
We find for $g\le 5$ a stratification of  depth $g-2$ of the moduli space of curves $\Mcal_g$  with the property that its strata are affine and the classes of their 
closures provide a $\QQ$-basis for the Chow ring of $\Mcal_g$. The first property confirms a conjecture of one of us. The way we establish the second property yields new (and simpler) proofs of theorems of Faber and Izadi which, taken together, amount to the statement that in this range the Chow ring is generated by the $\lambda$-class.
\end{abstract}
\maketitle

\section*{Introduction}

The starting point of the present note is a beautiful topological result by Harer, who in \cite{Harer:86} via a deep application of the theory of Strebel differentials showed  that the (coarse) moduli space of smooth algebraic curves of 
genus $g \ge 2$, $\Mcal_g(\CC)$,  has the homotopy type of a complex of dimension $4g-5$. A decade ago, the second author asked with 
Hain (\cite{HaiLoo:98}, Problem 6.5) whether there is a Lefschetz type of 
proof of this fact.  That would for instance  be accomplished by 
finding a stratification of $\Mcal_g$ with all strata affine subvarieties 
of codimension $\le g-1$. One reason to make this desirable is that the relatively  easy inductive algebraic approach developed  by Arbarello-Cornalba in \cite{ArbCor:98} to the rational cohomology of the Deligne--Mumford compactification $\overline{\mathcal{M}}_g(\CC)$ of $\Mcal_g(\CC)$ has Harer's theorem as the only nontrivial input from geometric topology.  But as explained in \cite{vakil:04}, the existence of such a stratification has many other interesting cohomological consequences as well.  

The classical candidate for an affine stratification of $\Mcal_g$
was the flag of subvarieties  defined by Arbarello  in terms of Weierstrass points in \cite{Arb:74} and \cite{Arb:78}. Its smallest stratum is the hyperelliptic locus and hence  affine
and  Theorem~1 of \cite{Fon:01} shows that the same is true for the open
stratum. But the intermediate strata not always 
seem to be affine (although they are easily seen to be quasi-affine).
Here we attempt to obtain an affine stratification by looking at properties of the linear system of  quadrics passing through the canonical image of a smooth projective genus $g$ curve. This also singles out the hyperelliptic locus as the smallest stratum, for if the curve is nonhyperelliptic, then  this linear system is of dimension $\binom{g-2}{2}$, whereas in the  hyperelliptic case its dimension is higher. But otherwise our approach yields strata that in general  differ from those of Arbarello. For instance, 
the complement of our codimension 2 stratum is for $g=4$ resp.\ $g=5$ the  locus of curves whose canonical model lies on a smooth quadric  resp.\  is an intersection of quadrics (or equivalently, that are not trigonal).
Although we do not offer a uniform receipe, we succeed in obtaining an affine stratification  of $\Mcal_g$ of depth $g-2$ for $g\le 5$ (we give in fact two such for $g=5$), thus answering  Problem 6.5 in \cite{HaiLoo:98} affirmatively in this range. 

The stratifications that we construct have one other nice property that, we feel, 
justifies calling them \emph{perfect}, namely that the strata have trivial Chow ring and
that their closures yield a $\QQ$-basis for $A^\pt(\Mcal_g)$.  In order to be precise, let us recall the structure of $A^\pt(\Mcal_g)$ in this range. If $f:\Ccal\to S$  is a smooth projective  family
of curves of genus $g\ge 2$, then its  \emph{$i$-th tautological class} $\kappa_i=\kappa_i(\Ccal/S)$ is $f_*( c_1(\omega_{\Ccal/S})^{i+1})\in A^i(S)$ (here $A^\pt (S)$ stands for the Chow ring of $S$ with  rational coefficients). According to Mumford \cite{Mum:83}, the subalgebra  of $A^\pt (S)$ generated by these classes (denoted by $\Rcal^\pt (\Ccal /S)$ or more often simply by $\Rcal^\pt (S)$)  contains the $i$th Chern class $\lambda_i=\lambda_i(\Ccal/S)$ of the Hodge bundle $f_*\omega_{\Ccal/S}$ with $\lambda_1$ being a non-zero multiple of $\kappa_1$.   A theorem of one of us \cite{Loo:95} implies that  $\lambda_1^{g-1}= 0$ always. On the other hand, 
Faber \cite{Faber:99} proved that for the universal family over the moduli stack $\Mcal_g$, $\lambda_1^{g-2}\not= 0$,  so that we have an injection
\[
\QQ [\lambda_1]/(\lambda_1^{g-1})\to \Rcal^\pt(\Mcal_g)\subset A^\pt(\Mcal_g).
\]
The composite  map $\QQ [\lambda_1]/(\lambda_1^{g-1})\to A^\pt(\Mcal_g)$ is for $g\le 5$ also onto; this is due to Faber \cite{Faber:90}  for $g=3,4$ and  to Izadi \cite{izadi} for $g=5$.  We find  in every codimension $k\le  g-2$ a single stratum and so the closure of that stratum must represent a non-zero multiple of $\lambda_1^k$. We  prove this  property independent of the two works just cited and thus obtain a new (and considerably simpler) proof that $A^\pt(\Mcal_g)=\QQ [\lambda_1]/(\lambda_1^{g-1})$ in this range. 
\\

In the first section we collect or  prove some general
results pertaining to $\Mcal_g$. As is well-known, a canonical curve possessing an effective even theta characteristic is contained in a rank 3 quadric. This is a codimension one phenomenon and so this suggests to take as our open stratum the locus of curves having no such theta characteristic. We indeed establish that  for $g \ge 3$ this open subset is affine and can be used as the open stratum for $g\le 5$.  The cases $g=3,4,5$ each have their own section. (The case $g=2$ is trivial: $\Mcal_2$ is affine and has trivial Chow ring.)
\\

We work here over a fixed  algebraically closed field $k$. We take it
to be of characteristic zero since some of the results we invoke are (still) only available
in that characteristic (the computations themselves only  seem  to exclude  characteristic 2 and 3).
\\

\emph{Notation.} If $C$ is a connected nonsingular projective curve, we write $V_C$ for $H^0(\omega_C)^*$ and denote by $f:C\to\PP V_C$ the  canonical map. We recall that when $C$ is hyperelliptic of genus $\ge 2$, $f(C)$ is a smooth rational curve of degree $g-1$ in the $(g-1)$-dimensional $\PP V_C$ and identifies the orbits of the hyperelliptic involution. Otherwise $f$ is an embedding with image a curve of degree $2g-2$ and we then identify $C$ with $f(C)$.

We also follow the custom to write $\lambda$ for $\lambda_1$. 
\\

We shall need the following basic result of geometric invariant theory (which is essentially due to Hilbert).

\begin{lemma}\label{quotient}
Let $U$ be an affine variety on which a reductive group $G$ 
operates.
If every $G$-orbit is closed (a condition that is obviously fulfilled if these orbits are fibers
of a morphism of varieties with domain $U$), then the orbit space $G\bs U$ is in a natural manner an affine variety. If in addition $G$ acts with finite stabilizers, then
$\dim (G\bs U)=\dim U-\dim G$. 
\end{lemma}
\begin{proof}
It is a classical theorem that,  $G$ being reductive,  $k[U]^G$ is noetherian
and separates the closed orbits. So the assumption that every orbit is closed implies that a closed point of $\spec (k[U]^G)$ corresponds to a $G$-orbit
in  $U$ and vice versa.  If $G$ also acts with finite stabilizers, then the  fibers of  
$U\to G\bs U$ will have the same dimension as $G$ and the last assertion easily follows.
\end{proof}

As a modest application of this lemma we recover  the well-known fact that the hyperelliptic locus $\Hcal_g\subset\Mcal_g$ is affine:
we may identify $\Hcal_g$ with $\Scal_{2g+2}\bs\Mcal_{0,2g+2}$ and 
$\Mcal_{0,2g+2}$ can be identified with an open affine subset of affine $(2g+1)$-space.
It is also well-known that $\Hcal_g$ is an irreducible and closed suborbifold of
$\Mcal_g$ of codimension $g-2$ and has trivial Chow ring. In particular, 
$\Mcal_2=\Hcal_2$ has these properties.

\section{Some genus independent results}

\subsection*{The Thetanull divisor} For $g\ge 3$, the locus in $\Mcal _g$ parametrizing curves with an effective even theta characteristic is an irreducible  hypersurface in $\Mcal_g$ that we shall refer to as the \emph{Thetanull} divisor and denote by 
$\Mcal'_g$.

\begin{proposition}\label{prop:top}
The  Thetanull divisor $\Mcal'_g$ is  irreducible and its complement in $\Mcal_g$ is affine for every 
$g \ge 4$.
\end{proposition} 
\begin{proof}
The irreducibility is due to Teixidor i Bigas \cite{Tei:88}, Proposition~(2.4). For the other property
we argue as for the proof of \cite{Fon:01}, Theorem~1. Indeed, if we 
recall that the rational Picard group of $\overline{\Mcal}_g$ is generated by the Hodge class 
$\lambda$ and  the classes $\delta_i$ of the irreducible components  
$\Delta_i$ ($i=0,\dots ,[\frac{1}{2}g]$) of the boundary  divisor $\partial \overline{\Mcal}_g$, then  according to
Teixidor i Bigas (\emph{op.~cit.}, Proposition~(3.1)), the class of  the closure $\overline{\Mcal}'_g$ 
of the Thetanull divisor in $\Mcal_g$ is $a\lambda - \sum_i b_i \delta_i$ with $a > 11$ and $b_i > 1$ 
for every $i$. As a consequence, the divisor
\[
D := \overline{\Mcal}'_g + \sum_i (b_i - 1) \Delta_i 
\]
is effective of class  $a \lambda - \delta$ and has support 
$\overline{\Mcal}'_g \cup \partial \overline{\Mcal}_g$. Since $a > 11$, we deduce from the 
Cornalba--Harris criterion  (\cite{CorHar:88}, Theorem 1.3) that $D$ is ample on 
$\overline{\Mcal}_g$. Hence the complement  of its support (which is also the complement 
of $\Mcal'_g$ in $\Mcal_g$) is affine.
\end{proof}

We shall find that  this proposition  also holds for $g=3$
(as is well-known, $\Mcal'_3=\Hcal_3$).

\subsection*{Half-canonical pencils and quadrics}

We  need the following lemma. 

\begin{lemma}\label{lemma:eventheta}
Let $C$ be a nonsingular projective curve of genus $g\ge 3$. 
Any effective even theta characteristic determines a rank 3 quadric in $\PP V$ that contains $f(C)$. Conversely, a rank 3 quadric $Q\subset \PP V$ whose smooth
part  contains $f(C)$ defines a half-canonical pencil on $C$.
\end{lemma}
\begin{proof}
Let $L,\phi:L^2\cong \omega_C$ be an effective even theta characteristic.
If $\alpha,\beta\in H^0(C,L)$ are linearly independent, then $u:=\phi(\alpha^2)$, $v:=\phi(\beta^2)$ and $w:=\phi(\alpha\beta)$ are elements of $V_C^*$ that satisfy $uv=w^2$.
The quadric defined by this equation evidently contains $f(C)$.

Now let  $Q\subset \PP V_C$ be a rank 3 quadric that contains $f(C)$. So its singular set
$Q_\sing$ is a codimension 3 linear subspace of $\PP V_C$. The dual variety 
$\check{Q}$ of $Q$ is a conic in the  projective plane in $\PP V^*_C$ dual to $Q_\sing$. Any $p\in \check{Q}$ defines a hyperplane in $\PP V_C$ that meets $Q$ in a codimension 2 linear subspace in $\PP V_C$  which  contains $Q_\sing$. Such a hyperplane  cuts out on $C$ a $2$-divisible divisor, i.e., one of the form 
$2D_p$. It is clear that  $D_p$ half-canonical  and moves in a pencil. 
\end{proof}

\subsection*{The trigonal locus} We denote by   $\Tcal_g\subset \Mcal_g$ the locus which parameterizes the genus $g$ curves admitting a $g^1_3$.
It is well-known (and easy to prove) that $\Tcal_g=\Mcal_g$ for $g=2,3,4$.
The following assertion is also well-known, but we give a 
proof for completeness.

\begin{proposition}\label{prop:trig}
For $g\ge 5$, the trigonal locus $\Tcal_g\subset \Mcal_g$ is a closed irreducible subvariety of 
$\Mcal_g$ of codimension $g-4$. It contains the  hyperelliptic locus $\Hcal_g$ (which is a subvariety of codimension $g-2$).
\end{proposition}
\begin{proof}
Denote by $\tilde\Tcal_g$ the moduli space of pairs $(C,P)$, where $C$ is a  nonsingular projective genus  curve of genus $g$ and $P$ is a $g^3_1$ on $C$. 
We first show that the forgetful map $\tilde\Tcal_g\to \Mcal_g$ is  proper. If
 $\Ccal_g\to \Mcal_g$ is  the universal genus $g$ curve  
and $\pic^3(\Ccal_g /\Mcal_g)$ the degree $3$ component of its Picard variety, then
we have an obvious morphism  $\tilde\Tcal_g\to \pic^3(\Ccal_g /\Mcal_g)$. The fiber over  a point $[(C,\ell )]$ of $\pic^3(\Ccal_g /\Mcal_g)$ is the Grassmannian of $2$-planes in $H^0(C,\ell)$. From this it is easily seen that  $\tilde\Tcal_g\to \pic^3(\Ccal_g /\Mcal_g)$ is proper. Hence so is $\tilde\Tcal_g\to \Mcal_g$.

The preimage of the hyperelliptic locus $\Hcal_g$ in $\tilde\Tcal_g$ can be identified with the universal hyperelliptic curve $\Ccal_{\Hcal_g}$: if $(C,P)$ represents a point of  $\tilde\Tcal_g$ and $P$ has a fixed point, then $C$ is hyperelliptic and the moving part of $P$ is the hyperelliptic pencil. Conversely, a hyperelliptic curve and a point on it determine a point of $\tilde\Tcal_g$ and so $\tilde\Tcal_g$ contains the universal hyperelliptic curve. 

We next show that a smooth projective curve $C$ of genus $g\ge 5$ cannot have two distinct $g^3_1$'s. Suppose it did so that we have a morphism 
$\phi: C\to \PP^1\to \PP^1$ of bidegree $(3,3)$. Since a  reduced curve in 
$\PP^1\times\PP^1$ of bidegree $(a,b)$ has for positive $a,b$ (by the adjunction formula) arithmetic genus resp.\  $(a-1)(b-1)$, $\phi$ cannot be of degree one onto its image.  It is then easy to see that  both pencils have a common fixed point.
But then they are hyperelliptic and coincide. This contradicts our assumption.

So $\tilde\Tcal_g\setminus\Ccal_{\Hcal_g}$ maps bijectivily  onto 
$\Tcal_g\setminus\Hcal_g$. A standard deformation argument shows that this is \'etale, hence it is an isomorphism.

Any element of $\Tcal_g\setminus\Hcal_g$ is represented
by a triple cover $C\to \PP^1$ that is unique up to an automorphism of $\PP^1$. The discriminant of this cover has degree  $2g+4$. From this we see that $\Tcal_g$
has dimension $2g+4-3=2g+1$.

According to Fulton \cite{Fulton:69}, the
moduli space of smooth connected triple coverings $ \PP^1$ of genus $g$  is irreducible and so $\Tcal_g$ is irreducible as well.
\end{proof}

\section{The case $g=3$} 
In the theorem below the second part is due to Faber \cite{Faber:90}. The proof given here
is not only simpler, but also has the virtue of yielding the first part (which is probably known).

\begin{theorem}
An affine stratification of $\Mcal_3$ of depth 1 is given by its hyperelliptic locus (which is
also its Thetanull divisor): $\Mcal_3\setminus\Hcal_3$ and $\Hcal_3$ are affine and irreducible.
Moreover the strata $\Mcal_3\setminus\Hcal_3$ and $\Hcal_3$ have the  Chow ring of a point  and so $A^\pt(\Mcal_3)$ is $\QQ$-spanned by  the classes of $\Mcal_3$ and $\Hcal_3$.
\end{theorem}
\begin{proof}
We know  that $\Hcal_3$ is affine and that $A^\pt(\Hcal_3)$ is spanned by $[\Hcal_3]$. So we must prove  that $\Mcal_3\setminus\Hcal$ is affine and $A^\pt(\Mcal_3\setminus\Hcal)$
is spanned by $[\Mcal_3-\Hcal]$.

Let $C$ be a nonhyperelliptic genus $3$ curve and let  $p_1,p_2\in C$ be in the support
of two different odd theta characteristics of $C$.
We shall use these theta characteristics to find normal coordinates in $\PP V_C$.
If we think of $C$ as canonically embedded as a quartic curve  in the projective plane $\PP V_C$, then
the projective tangent line $l_i\subset \PP V_C$ to $C$ at $p_i$ is a bitangent of $C$ and  $l_1\not=l_2$.
Denote the point of intersection of these  bitangents by $o$  and by $l_\infty$ the line through $p_1$ and $p_2$. It is clear that $o\notin C$. For $i=1,2$, the line $l_i$ meets $C$ in $p_i$ and a point $q_i$ (that may equal $p_i$) and $l_\infty$ meets
$C$ in $p_1,p_2$ and two other points, $r_1,r_2$, say (which may coincide but differ from $p_1$ and $p_2$). 
Choose coordinates
$[Z_0:Z_1:Z_2]$ such that $o=[1:0:0]$, $p_1=[0:1:0]$, $p_2=[0:0:1]$ and
$\{ r_1,r_2\}=\{[0:1:k],[0:k:1]\}$ for some $k\not= 0$. In terms of the affine coordinates 
$(z_1,z_2)=(Z_1/Z_0,Z_2/Z_0)$, an equation $f$ for $C$ has the form
\begin{multline*}
f(z_1,z_2)=1- 2a_1z_1 -2a_2z_2+a_1^2 z_1^2+a^2_2z_2^2 +bz_1z_2+\\
+(c_1z_1+c_2z_2)z_1z_2+ d(z_1^2+z_2^2)z_1z_2+ez_1^2z_2^2,
\end{multline*}
where $q_1=[a_1:0:1]$ and $q_2=[a_2:1:0]$ and $[1:k]$ and $[k:1]$ are solutions of $d(z_1^2+z_2^2)+ez_1z_2=0$.
Hence $C$ is given by  $(a_1,a_2,b,c_1,c_2,d,e)\in \AA^7$.
The affine coordinate system $(z_1,z_2)$ is unique up to the action of $\{ \pm 1\}\times \GG_m$ defined
by $(\varepsilon,t)\cdot (z_1,z_2)= (t^{-1}z_1,\varepsilon t^{-1} z_2)$ (the $\{ \pm 1\}$-action is needed because we did not order the pair $r_1,r_2$). The corresponding action on $\AA^7$ is given by
\[
(\varepsilon,t)\cdot (a_1,a_2,b,c_1,c_2,d,e)= 
(t^2a_1,t^2a_2,,\varepsilon t^2b,\varepsilon t^3c_1,t^3c_2,\varepsilon t^4d,t^4e).
\]
We use the discriminant to eliminate the $\GG_m$-action: 
For a polynomial $F$ homogeneous of degree $4$ in $Z_0,Z_1,Z_2$ we have defined 
its discriminant $\Delta(F)\in k$ (of degree $2^3 \cdot 28$). 
It is non-zero if and only if $F$ defines a nonsingular
quartic and it is a semi-invariant for the action of $\GL (3)$. So
its restriction to $\AA^7$, $\delta\in k[a_1,\dots ,e]$,  will
be a semi-invariant for the group $\{ \pm 1\}\times \GG_m$.
Denote by $\PP$  the weighted projective space $\GG_m\bs (\AA^7\setminus\{ 0\})$
and by $H\subset\PP$ the hypersurface defined by $\delta =0$.  Since  the morphism
$\PP \setminus H\to \Mcal_3\setminus \Hcal_3$ is finite, 
Lemma \ref{lemma:proj} below implies that $\Mcal_3\setminus \Hcal_3$ is affine  and $A^\pt(\Mcal_3\setminus \Hcal_3)=\QQ$.
\end{proof}

\begin{lemma}\label{lemma:proj}
Let $\PP$ be a weighted projective space (i.e., the proj of a graded polynomial algebra
whose generators have positive weights) and let $H\subset \PP$ be a hypersurface. Then $\PP\setminus H$ (and hence any variety admitting a finite cover isomorphic to $\PP\setminus H$) is affine and has the Chow ring of a point.
\end{lemma}
\begin{proof}
If $\PP=\proj k[z_0,\dots ,z_n]$ with $z_i$ having weight $d_i>0$, then 
the substitution $z_i=w_i^{d_i}$ defines a finite morphism  
$\PP^n\to \PP$. So without loss of generality we may assume that $\PP=\PP^n$.
Clearly, $\PP^n\setminus H$ is affine. Since $A^\pt (\PP^n)$ is generated by the hyperplane class, it is also generated by $[H]$ (which is a multiple of that class).
It follows that $A^\pt (\PP^n\setminus H)=\QQ$.
\end{proof}

Since we have an injective homomorphism   
$\QQ[\lambda]/(\lambda)^2\to A^\pt (\Mcal_3)$, we find:

\begin{corollary}
The map $\QQ[\lambda]/(\lambda)^2\to A^\pt (\Mcal_3)$ is an isomorphism and the
class $[\Hcal_3]\in A^1(\Mcal_3)$ is non-zero and a multiple of $\lambda$.
\end{corollary}

\section{The case $g=4$} 
According Lemma  \ref{lemma:eventheta} a smooth projective curve of genus 4 
possesses a half-canonical pencil if and only its the canonical model of $C$ lies on a singular conic. 
The corresponding locus $\Mcal'_4$ in $\Mcal_4$ is  the Thetanull divisor which by Proposition  
\ref{prop:trig} contains the hyperelliptic locus $\Hcal_4$.

\begin{theorem}
A stratification of $\Mcal_4$ of depth 2  is given by the Thetanull divisor and the hyperelliptic locus:
$\Mcal_4\supset \Mcal'_4\supset\Hcal_4$. This is a filtration by subvarieties whose successive 
differences are affine and have  the Chow ring of a point (so that $\Acal^\pt(\Mcal_4)$ is spanned by the classes
of $\Mcal_4, \Mcal'_4,\Hcal_4$).
\end{theorem}

Since we already know that $\Hcal_4$ is irreducible,  affine and has the Chow ring of a point, the theorem follows from
the two statements below.

\begin{proof}[Proof that $\Mcal_4\setminus \Mcal'_4$ is affine and has the Chow ring of a point]
Let be gi\-ven a projective curve $C$ of genus $4$ that represents a point of 
$\Mcal_4\setminus\Mcal'_4$, i.e., has no effective even theta characteristics. 
Then $C$ is in $\PP V_C$ the complete intersection of a nonsigular quadric surface $Q$ 
and a cubic surface. We recall that such a quadric  has two rulings (it is isomorphic to 
$\PP^1\times\PP^1$). We shift our attention to $Q$ and regard
$C$ as a divisor  on $Q$.

We choose three items attached to $C$: an ordering of the two rulings of $Q$, 
an effective degree $3$ divisor $D$ on $C$  that is an odd  theta characteristic and
a point $p$ in the support of $D$. Then  $D$ determines a plane section $k\subset Q$  with 
the property that within $Q$, $k$ meets $C$ in  $2D$. 
Let $m$ be the line  in the second ruling of $Q$ through $p$. This line and $C$ meet inside $Q$
at $p$ simply and hence they meet in two other (possibly coinciding) points  $\not= p$ as well. 
Denote by $l$ the line in the first ruling through the barycenter in $m\setminus\{ 0\}$ of these two  points. 
Now choose an affine
coordinate system $(x,y)$ for $Q$ such that the lines $k,l,m$ are defined by resp.\  $x=y$, $y=0$, $x=\infty$.
This coordinate system is unique up to the $\GG_m$-action defined by 
$t\cdot (x,y)= (t^{-1}x,t^{-1}y)$. If
\[
f(x,y)=\sum_{i,j=0}^3 a_{ij}x^iy^j
\]
is an equation for $C$, then $t\in\GG_m$ acts on $a_{ij}$ as multiplication by 
$t^{i+j}$. The conditions imposed  imply that $f(x,0)$ is of the form
$x(u+vx^2)$, in other words, $a_{0,0}=a_{2,0}=0$,  and that 
$f(x,x)$ is of the form $x^2(a+bx^2)^2$, in other words,
$\sum_{i+j=k} a_{i+j}$ equals $0$ for $k=1,3,5$, and $a^2, 2ab,b^2$ for $k=2,4,6$.
We may therefore parameterize such equations by affine $10$-space:
\[
\Big((a_{ij})_{i=1,2; j=1,2,3}, a_{1,0}, a_{3,0}, a,b\Big)\in\AA^{10}
\]
But beware: $(a^2,2ab,b^2)$ is read off from the equation, not $(a,b)$, and so
$(a,b)$ is defined up to a common sign. We give $a$ resp.\ $b$ weight $1$ resp.\ $3$
so that $\GG_m$ now acts on $\AA^{10}$ with positive weights.
>From here on we  proceed as in the case $g=3$: the  $\GG_m$-action on $\AA^{10}$ defines a weighted projective space $\PP$ of dimension $9$, 
the discriminant $\delta\in k[\AA^{10}]$ determines a hypersurface $H\subset \PP$. The morphism $\PP\setminus H\to \Mcal_4\setminus \Mcal'_4$ is finite.  It then follows from Lemma \ref{lemma:proj}  that  $\Mcal_4\setminus \Mcal'_4$ is affine and has trivial Chow ring.
\end{proof}

\begin{proof}[Proof that $\Mcal'_4\setminus \Hcal_4$ is affine and has the Chow ring of a point]
Let be given a projective curve $C$ that represents an element of $\Mcal'_4\setminus \Hcal_4$.
Then the canonical model of $C$  is complete intersection of
a quadric cone $Q$ and a cubic surface.  Denote the vertex of $Q$ by $o$. Since $C$ is smooth, 
$o\notin C$. 
We attach to $C$ both an effective degree $3$ divisor $D$ 
that is an odd theta characteristic and a point $p$ in the support of $D$. Since $2D$ is 
canonical, it is the  section of  a plane $H\subset \PP V_C$. This plane will not pass through $o$: 
otherwise $H$ would be tangent to $Q$, hence $D$ would lie in the half-canonical pencil defining $Q$, 
and this would contradict our assumption that $D$ is odd. Since the line $l$ through $o$ and $p$ is not 
contained in $H$, $l$ and  $C$, viewed
as curves on $Q$, meet simply in $p$. In particular, $l-\{ o,p\}$ meets $C$ in two points 
$p_1,p_2$ (which may  coincide). We put $Q_\infty:=Q\cap H$.

Now $C$ is defined by a regular function $f$ on $Q\setminus Q_\infty$ which has a pole of order three along 
$Q_\infty$. We normalize $f$  by requiring that $f(o)=1$.
Regard $\PP V_C\setminus H$ as a vector space with origin $o$ so that we can write
$f=f_3+f_2+f_1+1$ with $f_i$ homogeneous of degree $i$.
The quadratic part $f_2$  defines an effective degree $4$ divisor $E$ on $Q_\infty$. 
It follows from the preceding that the restriction of $f$ to $l\setminus \{ p\}$ must be of exact degree two 
(defining the divisor $(p_1)+(p_2)$) and so $p$ is not in the support of $E$. Regard $E$ as a divisor on  
$Q_\infty \setminus \{ p\}$ (which is isomorphic to an affine line) and denote by $m$ the line through 
the barycenter of $E$ in $Q_\infty \setminus \{ p\}$ and $o$.
We now  have  coordinate hyperplanes in $\PP V_C\setminus H$ spanned by
$l$ and $m$ and the two tangent planes of $Q$ containing $l$ and $m$. 

It is then easy to see that there exist  coordinates $(x,y,z)$ in $\PP V_C\setminus H$ such that 
$l$ is the $y$-axis, $m$ is the $x$-axis, $Q$ is given by $xy=z^2$ and 
$y(p_1)y(p_2)=1$. This system
of coordinates is unique up  to the action of the  group $\{\pm 1\}\times\GG_m$ given by
$(\eps ,t)\cdot (x,y,z)=(t^{-2}x,\eps y,t^{-1}z)$.
Since $Q_\infty\cdot C=2D$, $f_3$ can be given the form $x(ax+by+cz)^2$, where $ax+by+cz=0$ defines
in $H$ the line that is spanned by the degree $2$ divisor  $D-(p)$ (so the factor  $ax+by+cz$ is unique up 
to sign).
We represent $f_2$ as a polynomial without the monomial $z^2$. The condition imposed
on $f_2$ means that $f_2(1,t^2,t)$ is a quartic polynomial with vanishing $t^3$-coefficient.
This means that $f_2$ has vanishing $yz$-coefficient. The fact that $f(0,y,0)$ is quadratic  with 
constant term $1$ and the product of its roots equal to $1$ implies that the coefficient of $y^2$ is 
$1$. So we have now written
\[
f(x,y,z)=1+(a_2x+b_2y+c_2z)+(a_1x^2+b_1xy+y^2+c_1xz)+x(ax+by+cz)^2.
\]
These equations are parameterized by $(a,b,c, a_1,b_1,c_1,a_2,b_2,c_2)\in \AA^9$,
where we should bear in mind that $(-a,-b,-c; a_1,b_1,c_1,a_2,b_2,c_2)$ yields the same
equation. The $\{\pm 1\}\times\GG_m$-action lifts to
\begin{multline*}
(\eps ,t)\cdot (a,b,c; a_1,b_1,c_1,a_2,b_2,c_2)=\\
=(t^3a,\eps tb,t^2c; t^4a_1,\eps t^2b_1,t^3c_1,t^2a_2,\eps b_2,tc_2).
\end{multline*}
We proceed as before: the  $\GG_m$-action on $\AA^{9}$ defines a weighted projective space  $\PP$ of dimension $8$ and the discriminant $\delta\in k[\AA^{9}]$ determines a hypersurface  $H\subset \PP$ so that we have a finite morphism
$\PP\setminus H\to \Mcal'_4\setminus \Hcal_4$. It remains to apply Lemma 
\ref{lemma:proj}.
\end{proof}

Combining this with  the injectivity of  $\QQ[\lambda]/(\lambda)^3\to A^\pt (\Mcal_4)$
this yields:

\begin{corollary}
The map $\QQ[\lambda]/(\lambda)^3\to A^\pt (\Mcal_4)$ is an isomorphism and the
classes $[\Mcal'_4]$ resp.\ $[\Hcal_4]\in A^1(\Mcal_4)$ are non-zero multiples of 
$\lambda$ and $\lambda^2$ respectively.
\end{corollary}

In fact,  according to \cite{Tei:88}, Proposition~(3.1),  
$[\Mcal'_4]\in A^1(\Mcal_4)=34 \lambda$. 
On the other hand, we have by \cite{Mum:83}, \S 7,
$[\Hcal_4]=3 \kappa_2 - 15 \lambda^2$, and from the two formulas 
for the $\mathbb{Q}$-class of the hyperelliptic locus given in 
\cite{Ran:85}, Example 5.7 (a), we deduce 
the relation $2 \kappa_2=27 \lambda^2$ in the case $g=4$. 

\section{The case $g=5$} Let $C$ be a smooth projective curve $C$ of genus 5.
If $C$ is nonhyperelliptic, then the quadrics passing through its canonical model
form a net and either that net defines $C$ (so that $C$ is a complete intersection of three quadrics) or a rational scroll (for this and what follows, see \cite{stdonat}). 

In the last case, the scroll is a Hirzebruch surface of type $F_1$ ($=$ a projective plane blown up in a point) embedded by the linear  system $|e + 2f |$, where  $e\in\pic (F_1)$ is  the class of the exceptional curve $E\subset F_1$ and $f\in\pic (F_1)$  the class of the ruling, and we have $C\in | 3e+ 5f |$. So $C$ is trigonal: $C\cdot f=3$  and the ruling of the scroll yields a $g^1_3$ on $C$. It also follows that $C\cdot E=2$. So upon contracting $E$ in $F_1$ we get a projective plane curve
in which the image of $C$ becomes a quintic with either an ordinary double point
or a cusp (according to whether $C$ meets $E$ in two points or one). The converse is also true: a quintic plane curve $C'$ with a singular point of arithmetic genus one defines a genus 5 curve with a $g^1_3$ defined by the pencil of lines through that point. In the cusp case, 
our $g^1_3$ on $C$ is the moving part of a half-canonical pencil that has  the unique
point of intersection of $C$ and $E$ as its fixed point.  (According to Accola \cite{Accola:74}, any half-canonical pencil on a trigonal curve is necessarily of this form and so $C$ has no such pencils if $C'$ has an ordinary double point.) 

Recall from Proposition \ref{prop:trig} that the trigonal locus $ \Tcal_5\subset\Mcal_5$ contains the hyperelliptic locus.  Since a hyperelliptic curve 
possesses half-canonical pencils, we have in fact $\Hcal_5\subset \Tcal_5\cap\Mcal'_5$. The main result of this section is:

\begin{theorem}\label{thm:g=5}
An affine stratification of $\Mcal_5$ of depth 3 consists of the Thetanull divisor, the intersection of the trigonal locus with the Thetanull divisor and the hyperelliptic locus; 
another one is obtained by replacing the Thetanull divisor by the trigonal locus.
So if we put $\Tcal'_5:= \Tcal_5\cap \Mcal'_5$, then these are given by the following filtrations by closed irreducible subvarieties:
\begin{gather*}
\Mcal_5\supset \Mcal'_5\supset \Tcal'_5\supset \Hcal_5\text{ resp. }
\Mcal_5\supset \Tcal_5\supset \Tcal'_5\supset \Hcal_5.
\end{gather*}
Moreover, the successive differences $\Mcal_5\setminus \Tcal_5$,
$\Tcal_5\setminus\Tcal'_5$ and
$\Tcal'_5\setminus\Hcal_5$ have trivial Chow ring.
\end{theorem}

We know that $\Hcal_5$ is a closed irreducible subvariety of $\Mcal_5$ that is also affine. Following Proposition \ref{prop:top}, $\Mcal'_5$ is an irreducible hypersurface in $\Mcal_5$ whose complement is affine. We discuss the remaining strata separately.

The proof that the Chow ring of $\Mcal_5\setminus\Tcal_5$ is tautological, requires
a few general remarks about Chern classes attached to families of curves.
Recall that if $C$ is a nonhyperelliptic projective smooth curve, then the natural map 
$\sym^2 H^0(C,\omega_C)\to H^0(C,\omega^2_C)$ is surjective and its kernel defines the linear system of quadrics through $C$ in $\PP V_C$. 

\begin{lemma}\label{lemma:quadratickernel}
Let $f:\Ccal \to S$ be a smooth projective family of nonhyperelliptic  curves. Then we have a corresponding surjection of vector bundles
\[
\sym^2 f_*\omega_{\Ccal /S}\to f_*\omega^2_{\Ccal /S}
\]
whose kernel  has the property that  its Chern classes (taken in the Chow ring of $S$) are in the subalgebra generated by the tautological classes $\kappa_i(\Ccal/S)$.
\end{lemma}
\begin{proof} The $i$-th Chern class $\lambda_i$ of the Hodge bundle $f_*\omega_{\Ccal /S}$ is according to Mumford a universal expression in the $\kappa$-classes.
Since the Chern classes of a symmetric power of a vector bundle  are universally expressed in terms of those of the vector bundle, the same is true for 
$\sym^2 f_*\omega_{\Ccal /S}$. It therefore remains to show the corresponding property for $f_*\omega^2_{\Ccal /S}$. 

The proof of this is a minor variation of Mumford's. Since  the higher direct images
of $\omega^2_{\Ccal /S}$ vanish, the total Chern class of the bundle $f_*\omega^2_{\Ccal /S}$ can be computed by means of Grothendieck-Riemann-Roch: 
it tells us that this Chern class is the direct image on $S$ of the product of the relative Todd class of $f$ and the Chern character
of $\omega^2_{\Ccal /S}$.  The Todd class  is a polynomial in 
$c_1(\omega_{\Ccal /S})$  and so is  the Chern character of $\omega^2_{\Ccal /S}$
(for it equals $\exp (2c_1\omega_{\Ccal /S})$).
The desired property then follows. 
\end{proof}

\begin{proof}[Proof that $\Mcal_5\setminus \Tcal_5$ is affine and has trivial Chow ring] Let $V$ be  a fixed vector space of dimension $5$ and denote by $\sym^2 V^*$ the
$15$-dimensional space of quadratic forms on $V$. A net of quadrics in $\PP V$
is given by a $3$-plane in $\sym^2 V^*$. By the purity of the discriminant, the nets that fail to define a smooth complete intersection make up a hypersurface $H$ in the  Grassmannian $G_3(\sym^2 V^*)$ of $3$-planes in $\sym^2 V^*$.  The Picard group of  $G_3(\sym^2 V^*)$ is generated by an ample class and so the complement $U:= G_3(\sym^2 V^*)\setminus H$ is affine and $A^1(U)=0$. Since $U$ parametrizes canonically embedded genus 5 curves, we have a morphism $U\to \Mcal_5$ whose fibers are $\SL (V)$-orbits. It then follows from Lemma \ref{quotient} that $\SL (V)\bs U$ is affine. The image of $U$ in $\Mcal_5$ is clearly $\Mcal_5\setminus  \Tcal_5$. So $\Mcal_5\setminus  \Tcal_5$ is affine.

Next we prove the tautological nature of the Chow ring. Notice that $U$ carries a family of genus $5$ curves to which Lemma \ref{lemma:quadratickernel} applies. So $A^\pt (U)$ is generated by the tautological classes of this family. In other words, the  map
$U\to  \SL (V)\bs U\cong \Mcal_5\setminus \Tcal_5$
induces a surjection of $\Rcal^\pt (\Mcal_5\setminus \Tcal_5)\to A^\pt (U)$.

Since the group $\SL (V)$ acts on $U$ with finite stabilizers and has 
$\Mcal_5\setminus  \Tcal_5$ as orbit space, it follows from a theorem of Vistoli \cite{vistoli:87} that the kernel of the algebra homomorphism $A^\pt (\Mcal_5\setminus  \Tcal_5)\to A^\pt (U)$ is the $\Rcal^\pt (\Mcal_5\setminus \Tcal_5)$-submodule generated by the `Chern classes' of the almost $\SL (V)$-principal bundle $U\to \Mcal_5\setminus  \Tcal_5$. These Chern classes are readily identified as the Chern classes of the
Hodge bundle over $\Mcal_5\setminus  \Tcal_5$ and so the latter's Chow ring
is tautological as asserted. Since the first Chern class of a $\SL (V)$-bundle is trivial,
it also follows that $A^1(\Mcal_5\setminus  \Tcal_5)=0$, in particular $\lambda=0$ in $A^\pt(\Mcal_5\setminus  \Tcal_5)$.

Now in order to conclude we may argue as in \cite{izadi}, proof of Corollary 3.3 
on p. 293: by \cite{Mum:83} p. 309 and \cite{Faber:90} p. 447, 
$\Rcal^\pt(\Mcal_5\setminus  \Tcal_5)$ is generated by $\kappa_2$. 
On the other hand, from formula (7.7) in \cite{Mum:83} by keeping 
into account the vanishing of $\lambda$ it follows that the class of 
the locus of curves with a point $p$ with $h^0(3p) \ge 2$ is $81 \kappa_2$.
Since this locus is contained in $\Tcal_5$, its class vanishes in 
$A^\pt(\Mcal_5\setminus  \Tcal_5)$, hence our claim follows. (Alternatively, 
we can invoke  Faber's verification in   \cite{Faber:99} that
$\Rcal(\Mcal_5)$ is generated by $\lambda$. Since $\lambda|\Mcal_5\setminus  \Tcal_5=0$, it follows that $\Rcal^\pt(\Mcal_5\setminus  \Tcal_5)=\QQ$.)
\end{proof}

\begin{proof}[Proof that $\Tcal_5\setminus \Tcal'_5$ and $\Tcal'_5\setminus\Hcal_5$ are affine and have trivial Chow ring]
Let $C'$\\ be a quintic plane curve in a projective plane $P$ whose singular set  consists of a node or a cusp, say at $p\in P$. 
A generic  line $\ell$ in $P$ through $p$ meets $C'$ in three other points.  
Then it is easily verified that the barycenter in $\ell$ of these three points in the affine line $\ell\setminus \{ p\}$ 
traces, as $\ell$  varies, a cubic curve $K\subset P$ which has the same $3$-jet in $p$ as  $C'$. If $K$ is a nodal cubic  curve, then it has  has $3$ flex points, otherwise there is just one. Choose a flex point and denote by $l$ the  the tangent to $K\setminus \{ p\}$ at its unique flex point.
There exist affine coordinates  $(x,y)$ for $P\setminus l$ such that $K$ is given by 
$y^2+x^2+x^3$ resp.\ $y^2+x^3=0$. In the nodal  case these this system is unique up
to a sign: $(x,-y)$ works also, in the cuspidal case it is unique 
up to the $\{ \pm 1\}\times \GG_m$-action defined by
$(\eps ,t)\cdot (x,y)= (t^{-2}x,\eps t^{-3}y)$. Let us continue with the cuspidal case (the nodal case being easier). We then find  that $C'$ has an equation of the form
\[
x^3+y^2+\sum_{i+j=4,5} a_{ij} x^iy^j.
\]
The coefficient space is an $11$-dimensional affine space $\AA^{11}$. The group  
$\{ \pm 1\}\times \GG_m$ acts on this space componentwise by letting
$(\eps ,t)$ send $a_{ij}$ to  $\eps^jt^{2i+3j}a_{ij}$. The $\GG_m$-action  defines a weighted projective space $\PP$ and the discriminant yields a hypersurface $H\subset \PP$ such that we get a finite morphism $\PP\setminus H\to \Tcal'_5\setminus\Hcal_5$. The assertion then follows from Lemma \ref{lemma:proj}.
\end{proof}

\begin{corollary}
The map $\QQ[\lambda]/(\lambda^4)\to A^\pt(\Mcal_5)$ is an isomorphism and
for each of the two filtrations of Theorem \ref{thm:g=5}, the classes of its members define
a $\QQ$-basis for $A^\pt (\Mcal_5)$. 
\end{corollary}
\begin{proof} We know that the map  $\QQ[\lambda]/(\lambda^4)\to A^\pt(\Mcal_5)$
is injective.
Since the filtration $\Mcal_5\supset\Tcal_5\supset\Tcal'_5\supset\Hcal_5$ 
consists of varieties with the property that the restriction of $\lambda$ to each
successive difference is zero and are irreducible (except the last one), we have
that $\lambda^k$ is proportional to the class of the $k$th member of this filtration.
The theorem of Teixidor i Bigas \cite{Tei:88}  cited earlier implies that this remains true
if we replace $\Tcal_5$ by $\Mcal'_5$.
\end{proof}

\begin{question}
Since for $g\ge 5$, $\Tcal_g$ contains  $\Hcal_g$ as a closed affine subvariety of codimension $2$, 
we wonder whether  for such $g$ there exists an affine stratification of  $\Tcal_g$  of depth $2$. 
The answer may very well be \emph{no}, but the preceding theorem shows  it is \emph{yes} when $g=5$.
\end{question}

\end{document}